\documentclass{article}

\usepackage{arxiv}
\usepackage{amsmath}
\usepackage[utf8]{inputenc} % allow utf-8 input
\usepackage[T1]{fontenc}    % use 8-bit T1 fonts
\usepackage{hyperref}       % hyperlinks
\usepackage{url}            % simple URL typesetting
\usepackage{booktabs}       % professional-quality tables
\usepackage{amsfonts}       % blackboard math symbols
\usepackage{nicefrac}       % compact symbols for 1/2, etc.
\usepackage{microtype}      % microtypography
\usepackage{lipsum}
\usepackage{multicol}
\usepackage[english]{babel}
\usepackage{amsthm}
\usepackage{graphicx}
\usepackage{subfig}

\usepackage[capitalise,nameinlink,noabbrev]{cleveref}
\crefname{equation}{}{}

\newtheorem{definition}{Definition}[section]
\newcommand{\BF}[1]{\textbf{#1}}
\newcommand{\bb}[1]{\mathbb{#1}}

\theoremstyle{definition}

\theoremstyle{remark}

\newtheorem{exmp}{Example}[section]

\title{Division and new multiplication between vectors}

\author{
 	J. Enrique H. Ramírez \\
 Department of Physics\\
 ABC Federal University\\
 Sao Paulo, Brazil\\
 \texttt{enrique.ramirez@ufabc.edu.br}
 \And
 E Roca O \\
 Department of Physics\\
 University of Oriente\\
 Santiago de Cuba, Cuba\\
 \texttt{eroca@cnt.uo.edu.cu}
 
}

\begin{document}
\maketitle

\begin{abstract}
The division between two vectors belonging to the same vector space is obtained by elementary procedures of vector algebra and is defined by a matrix. This representation is obtained for two and three dimensional vector spaces. A new vector multiplication is defined and an the inverse vector. Through this multiplication we can obtain the division previously defined.The meaning of vector division and multiplication are analyzed.
\end{abstract}

% keywords can be removed
\keywords{Vector Division \and Vector Multiplication \and 
 Inverse Vector.}

\begin{multicols*}{2}
\section{Introduction}

The division between vectors in the framework of vector algebra is a subject not widely treated. Few studies have been carried out about this subject and those that have been carried out conclude that it is impossible to define this operation. In $1966$, $\text{Kenneth O. May}$ wrote an article,\cite{may1966impossibility}, in which he states that it is impossible to divide two vectors in a three-dimensional space.  This impossibility is essentially due to the ways that have been used to define division between vectors, which have not achieved this goal because they have tried to obtain it by one of the existing products between vectors, either by the dot product or by the vector product, these products not being proper multiplications like those defined between numbers. These ideas lead to endless solutions. In this work we will show that it is possible to define the division between vectors within the framework of vector algebra by means of a geometric representation.

	As is well known, in the field $\BF{K}$, of real numbers, $\bb{R}$, or of complexes, $\bb{C}$, the 	division operation can be entered as the inverse operation of multiplication, i.e., $a /b$ is true c such that $b c = a$; provided that $b \neq 0$, $c$ exists and is unique. If $b = 0$, $c$ exists, then and only if $a = 0$, and in this case $c$ is any number in $\BF{K}$. Since in the case of vectors there is no multiplication operation by which to perform the
	inverse operation, it is common to define the quotient between vectors as a particular tensor. Often the quotient between two quantities is of a different type and more complicated than this, the quotient of two integers is generally not an integer but a rational number. Similarly, the quotient of two vectors is not a vector, as is well known, but a quantity of a new type: a tensor \cite{goldstein2018mecanica}; but this means that the division of vectors can be exploited only with knowledge of the tensor algebra. In \cite{hestenes1971vectors}, division of vectors is treated in the context of multivector algebra, which is not treated in the same way in basic courses. In this algebra multiplication, and the inverse of a multivector, defines the division between multivectors which are vectors in a special case. 
	Our main goal in this paper is to use elementary procedures to divide 
     two vectors, $\BF{a}$ and $\BF{b}$, of a certain vector space $V$ of finite dimension with defined inner product. In what follows we will assume that $\BF{b} \neq \BF{0}$ unless otherwise stated, and 
	we denote by $ \vert\BF{a}\vert = a$, etc.

\section{Vector Division}

\begin{figure}[h]
	\centering
	\subfloat[Decomposition of the vector $\BF{a}$ according to the direction \\
	parallel and perpendicular to $\BF{b}$]{
		\label{fig1a}
		\hspace{1cm}
		\includegraphics[width=0.4\linewidth,height=0.26\textheight]{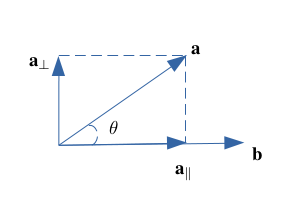}}
	\subfloat[The vector $\BF{b}$, its inverse ($\BF{b}^{-1}$), and the two possible vectors perpendicular to it ($\BF{b}_{\perp1}$, $\BF{b}_{\perp2}$) are shown. The inverse perpendicular vectors are shown.]{
		\label{fig1b}
		\includegraphics[width=0.4\linewidth,height=0.25\textheight]{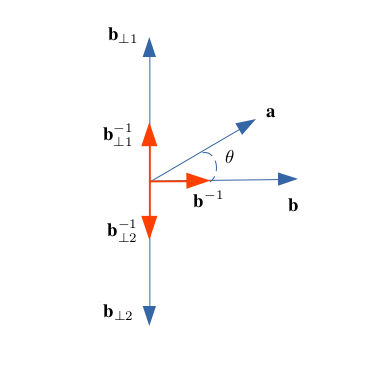}}
		\caption{The Fig. \ref{fig1a} show the decomposition of the vector $\BF{a}$ that allows defining the division between the vectors $\BF{a}$ and $\BF{b}$ and Fig. \ref{fig1b} the inverse and perpendicular vectors of the vector $\BF{b}$.}
	\label{vectoresdescomposition}
\end{figure}

Let $\BF{a}$ and $\BF{b}$ two vectors in finite n-dimensional space $V$, over the real numbers , with inner producto defined. 
We can divide vector $\BF{a}$ into two directions, one parallel  and the other perpendicular to vector the $\BF{b}$ as,

\begin{eqnarray}
\BF{a} = \BF{a}_{\|} +\BF{a}_{\perp},
\end{eqnarray}

thus, $\BF{a}_{\|} = \alpha\, \BF{b}$ and $\BF{a}_{\perp} = \beta\, \BF{b}_{\perp}$(see Fig. \ref{fig1a} ), where $\alpha$ and $\beta$ are two real proportionality coefficients to be determined, and $\BF{b}_{\perp}$  a perpendicular rotation of the vector $\BF{b}$. So we can write the above equation as,

\begin{eqnarray}\label{vectordecomposition1}
\BF{a} =  \alpha\, \BF{b} +\beta\, \BF{b}_{\perp}.
\end{eqnarray}

 Scalar multiplication by $\BF{b}$ and $\BF{b}_{\perp}$ let us to find the $\alpha$ and $\beta$ coefficients,
 
 \begin{eqnarray}
 \alpha &=& \frac{\BF{a}\cdot \BF{b}}{b^2},\nonumber\\
 \beta &=& \frac{\BF{a}\cdot \BF{b}_{\perp}}{b^2}.\nonumber
 \end{eqnarray}

We can write $\BF{b}_{\perp}$ as $\BF{b}_{\perp} = R\,\BF{b}$ , were $R$ is a perpendicular rotation matrix, so \ref{vectordecomposition1} can be writen as,

\begin{eqnarray}\label{vectordecomposition2}
\BF{a} =  \alpha\, I\, \BF{b} +\beta\, R\, \BF{b},
\end{eqnarray}
 
 where $I$ is the unit matrix. Factorization of $\BF{b}$ in \ref{vectordecomposition2} give,
 
 \begin{eqnarray}\label{vectordecomposition3}
 \BF{a} =  \left(\alpha\, I\ +\beta\, R \right)\,\BF{b},
 \end{eqnarray}

were we can define,

\begin{definition}
	The division between the vectors $\BF{a}$ and $\BF{b}$, which we denote as $\BF{a}/\BF{b}$, where $\BF{b} \neq \BF{0}$, is given by the matrix $E$,
	
	\begin{eqnarray}\label{defdivision}
	\frac{\BF{a}}{\BF{b}} &=& \alpha\, I\ +\beta\, R .
	\end{eqnarray}
\end{definition}

	If $V$ is the vector space $\mathcal{R}^2$, then

	\begin{eqnarray}	
E &= &	\left[ \begin{array}{clr}
	\alpha &  -\beta \\
	 \beta
	& \alpha 
	\end{array}\right]\\	
	&= &	\frac{a}{b}\left[ \begin{array}{clr}
	\cos\theta &  -\sin\theta \\
	 \sin\theta
	& \cos\theta
	\end{array}\right].	\nonumber
	\end{eqnarray}

	      Matrix $E$ rotates vector $\BF{b}$ and then adjusts its magnitude to the magnitude of vector $\BF{a}$. In the particular case that $\BF{a} = \BF{b}$, then $E$ is the usual rotation.
	      
	      In this case the determinant of $E$ (det$E$) e given by
	      
	   \begin{eqnarray}
	      det E &=& \alpha^2 + \beta^2\\ &=&\left(\frac{a}{b}\right)^2\nonumber
	   \end{eqnarray}
	      
	 	If $V$ is the vector space $\mathcal{R}^3$, we have
	 	     	\begin{eqnarray}	
	 	     E &= &	\left[ \begin{array}{clrr}
	 	     \alpha &  C\beta & - B\beta\\
	 	     -C\beta & \alpha & A\beta\\
	 	      B\beta& - A\beta& \alpha
	 	     	 	     \end{array}\right]\\
	 	     	 	     &= &\frac{a}{b}	\left[ \begin{array}{clrr}
	 	     	 	     \cos\theta &  C\sin\theta & - B\sin\theta\\
	 	     	 	     -C\sin\theta & \cos\theta & A\sin\theta\\
	 	     	 	     B\sin\theta& - A\sin\theta& \cos\theta
	 	     	 	     \end{array}\right],	\nonumber
	 	     	\end{eqnarray}
    where $A$, $B$ and $C$ are the components of the unit vector $\BF{n}$ given by,

\begin{eqnarray}
\BF{n} &=& \frac{\BF{a}\times\BF{b}}{\|\BF{a}\times \BF{b}\|}\\
&=& (A, B, C)\nonumber
\end{eqnarray}

Here,
 \begin{eqnarray}
det E &=& \alpha\left(\alpha^2 + \beta^2\right)\\ &=&\left(\frac{a}{b}\right)^3\cos\theta. \nonumber
\end{eqnarray}

Let's consider two examples:

\begin{exmp}\label{example1}
	Let \BF{a} = (3, -1) and \BF{b} = (2, 5) be two vectors in $\mathbb{R}^2$ find \BF{a}/\BF{b}.
	
	\begin{eqnarray*}
	\alpha &=& \frac{\BF{a}\cdot \BF{b}}{b^2}\\
	&=&\frac{1}{29}.		
	\end{eqnarray*}
\begin{eqnarray}
	\beta &=& \frac{\BF{a}\cdot\BF{b}_{\perp}}{b^2}\nonumber
\end{eqnarray}
There are two vectors perpendicular to $\BF{b}$, $\BF{b}_{1\perp} = (-5,2)$ and $\BF{b}_{2\perp} = (5,-2)$. The perpendicular rotation matrix $R$ will depend on which vector we choose. Let's choose $\BF{b}_{1\perp}$, then
\begin{eqnarray}
\beta &=& -\frac{17}{29}\nonumber
\end{eqnarray}
and

\begin{eqnarray*}
E &=& \frac{1}{29} I - \frac{17}{29} R\\
&=& 	\frac{1}{29}\left[ \begin{array}{clr}
	1 &  17 \\
	-17
	& 1
\end{array}\right].	\nonumber
\end{eqnarray*}

where $R$ is given by

\begin{eqnarray}
R  &=& 	\left[ \begin{array}{clr}
0 &  -1 \\
1
& 0
\end{array}\right]\nonumber
\end{eqnarray}
Note that if we choose $\BF{b}_{2\perp}$, the corresponding perpendicular rotation matrix would be the transpose of the previous one, on the other hand we can easily verify that  $\BF{a} = E\,\BF{b}$.\\

\end{exmp}

\begin{exmp}\label{example2}
	Let \BF{a} = (3, -1, 2) and \BF{b} = (2, 5, 1) be two vectors in $\mathbb{R}^3$ find $\BF{a/b}$.\\
	
In this case we find that $\alpha = 3/30$.   To find $\beta$, let's take the cross product of $\BF{a}$ and $\BF{b}$
\begin{eqnarray}
\BF{a}\times\BF{b} &=& -11 \BF{i} + \BF{j} +17 \BF{k}\nonumber
\end{eqnarray}

being
\begin{eqnarray}
\BF{n} &=& \frac{1}{\sqrt{411}}(-11,1,17)\nonumber
\end{eqnarray}
 
 In this way we find that $E$ is given by
 
 \begin{eqnarray}
 E  &= &	\frac{1}{30}\left[ \begin{array}{clrr}
 3 &  17 & -1\\
 -17 & 3 & -11\\
 1& 11& 3
 \end{array}\right].
 \end{eqnarray}
 
 Again we can verify that $\BF{a} = E\,\BF{b}$.
\end{exmp}
\subsection{Division Properties}

\begin{enumerate}
	\item [D1:]
	 $\BF{a}/\BF{a} = I$
	 \begin{proof}
	 	By  \cref{defdivision} we have that
	 	
	 	 $\BF{a}/\BF{a} = \BF{a}\cdot\BF{a}/a^2\,I + \BF{a}\cdot\BF{a}_{\perp}/a^2\,R$ but\\
	 	 $\BF{a}\cdot\BF{a}_{\perp} = 0$ and $\BF{a}\cdot\BF{a}/a^2 = 1$, then $\BF{a}/\BF{a} = I$.
	 \end{proof}
 \item [D2:]
 $\BF{0}/\BF{a} = O$
 
 \begin{proof}
 	By  \cref{defdivision} we have that
 	
 	$\BF{0}/\BF{a} = \BF{0}\cdot\BF{a}/a^2\,I + \BF{0}\cdot\BF{a}_{\perp}/a^2\,R$ then $\BF{a}/\BF{a} = O$.
 \end{proof}
\item[D3:]
$\BF{a}/\BF{b} = \left(\BF{b}/\BF{a}\right)^{-1}$
\begin{proof}
	We are going to prove this property by taking the product
	$(\BF{a}/\BF{a})(\BF{b}/\BF{a})$ and we will show that it is equal to the unitary matrix $I$.
	
	$(\BF{a}/\BF{b})(\BF{b}/\BF{a}) =(\BF{a}\cdot\BF{b}/b^2\,I+\BF{a}\cdot\BF{b}_{\perp}/b^2\,R)(\BF{a}\cdot\BF{b}/a^2\,I+\BF{a}\cdot\BF{b}_{\perp}/a^2\,\acute{R}) $\\
	
	Carrying out the indicated operation and considering that $R\,\acute{R} = I$ and $R = -\acute{R}$ we arrive at the expected result.
	
\end{proof}
\item[D4:]

$\BF{(a + b)}/\BF{c} = \BF{a}/\BF{c} + \BF{b}/\BF{c}$

\begin{proof}
	By definition,
	
	$\BF{(a + b)}/\BF{c} = ((\BF{a + b})\cdot\BF{c})/c^2\,I + ((\BF{a + b})\cdot\BF{c}_{\perp})/c^2\,R$\\
	
	Because scalar product follow distributive law we have that, 
	
	$\BF{(a + b)}/\BF{c} = \BF{a}\cdot\BF{c}/c^2\,I + \BF{ b}\cdot\BF{c}/c^2\, I + \BF{a}\cdot\BF{c}_{\perp}/c^2\,R + \BF{ b}\cdot\BF{c}_{\perp}/c^2\,R$
	
	Rearranging the terms we get that
	
	$\BF{(a + b)}/\BF{c} = \BF{a}/\BF{c} + \BF{b}/\BF{c}$, as we hope.

\end{proof}
 
\end{enumerate}

\section{Vector Multiplication}

The deficiency of thedot and cross products is due to the fact that through them it is not possible to define the inverse vector. For this reason, we will dedicate this section to defining a new multiplication between vectors, through which we can obtain the inverse vector, and that said multiplication be consistent with the division operation between vectors previously defined.\\

\begin{definition}[Vector Multiplication]\label{defmulti}
	Let  $V$ be a vector space with inner product defined and  $\BF{a}$ and $\BF{b}$ two vectors $\in$ $V$. The product between the vectors $\BF{a}$ and $\BF{b}$, which we denote as $\BF{a}\otimes\BF{b}$, is given by the matrix $E$,
	
	\begin{eqnarray}
	\BF{a}\otimes\BF{b} &=& E\\
	&=&\alpha\, I\ +\beta\, R\nonumber ,
	\end{eqnarray}

	where $I$ is the identity matrix, $R$  a perpendicular rotation matrix and $\alpha$ and $\beta$ two coefficients given by
	\begin{eqnarray*}
		\alpha &=&  \BF{a}\cdot \BF{b}\\
		 \beta &=&  \BF{a}\cdot \BF{b}_{\perp} = \BF{a}_{\perp}\cdot \BF{b}.
	\end{eqnarray*}

\end{definition}

The $\beta$ coefficient is closely related to the perpendicular rotation matrix $R$ and to the perpendicular vectors $\BF{b}_{\perp}$ or $\BF{a}_{\perp}$.

Now we are going to define the inverse vector by considering  \cref{defmulti} .\\

\begin{definition}\label{invvector}
	Let $\BF{a} \neq \BF{0} $  a vector in $V$. The inverse vector of $\BF{a}$, which we denote as $\BF{a}^{-1}$, is , from\,  \cref{defmulti} , such that
	
	\begin{eqnarray}
	\BF{a}\otimes\BF{a}^{-1} &=&  I.
	\end{eqnarray}
\end{definition}

From \cref{invvector} it follows that $\BF{a}^{-1} = \BF{a}/a^2$, which shows that the vectors $\BF{a}$ and $\BF{a}^{-1}$ are collinear. Now let's see the properties of multiplication.

For two and tree dimensional space we have that

\begin{eqnarray}
\BF{a}\otimes\BF{b} &=& a\,b 	\left[ \begin{array}{clr}
\cos\theta &  -\sin\theta \\
 \sin\theta
& \cos\theta 
\end{array}\right],
\end{eqnarray}

and

 	\begin{eqnarray}\label{multmatrix}	
\BF{a}\otimes\BF{b}&= &	\left[ \begin{array}{clrr}
\alpha &  C\beta & - B\beta\\
-C\beta & \alpha & A\beta\\
B\beta& - A\beta& \alpha
\end{array}\right]\\
&=& a\,b	\left[ \begin{array}{clrr}
\cos\theta &  C\sin\theta & - B\sin\theta\\
-C\sin\theta & \cos\theta & A\sin\theta\\
B\sin\theta& - A\sin\theta& \cos\theta
\end{array}\right]\nonumber
\end{eqnarray}

These matrices rotate the inverse vector $\BF{a}^{-1}$ ($\BF{b}^{-1}$) towards $\BF{b}$($\BF{a}$) and then adjust its magnitude converting it into the vector   $\BF{b}$ ( $\BF{a}$).
 
 From Eq.(\ref{multmatrix}) we can obtain the well-known vector algebra relation,
 
  \begin{eqnarray}
  (a\,b)^2 &=&\left(\BF{a}\cdot\BF{b}\right)^2 + \left(\|\BF{a}\times\BF{b}\|\right)^2.
 \end{eqnarray}
 
\subsection{Multiplication Properties}

\begin{enumerate}

	\item[M1:] $\BF{a}\otimes\BF{a}^-1{} = I$, multiplicative inverse.
	
	\begin{proof}
		We will show that there exists a multiplicative  inverse of $\BF{a}$ (which we denote as $\BF{a}^{-1}$) such that $\BF{a}\otimes\BF{a}^{-1} = I$\\
		
		Since $\alpha = 1$ that means that $a^{-1} = 1/a$. On the other hand  because $\beta = 0$ the vectors are collinear (that is, $\BF{a}^{-1} = k\,\BF{a}$), so we have that $k = 1/a^2$, in this way we find that $\BF{a}^{-1} = \BF{a}/a^2$.
		
	\end{proof}

\item[M2:]
	 $\BF{a}\otimes\BF{b} = \BF{b}\otimes\BF{a}$,  commutative law for multiplication.
	 \begin{proof}

	 	Let  $\BF{a}\otimes\BF{b} =\alpha\, I\ +\beta\, R,$ and  $\BF{b}\otimes\BF{a} =\acute{\alpha}\, I\ +\acute{\beta}\, \acute{R} $, with $\BF{a}\neq  \BF{0}$ and $\BF{b}\neq \BF{0}$ and $\theta$ the angle between the vectors $\BF{a}$ and $\BF{b}$. We must show that $\alpha = \acute{\alpha}$, $\beta = \acute{\beta}$ and $R = \acute{R}$.
	 	
	 	By definition $\alpha = \BF{a}\cdot\BF{b}$ and $\acute{\alpha} = \BF{a}\cdot\BF{b}$, so $\alpha = \acute{\alpha}$. On the other hand,  $\beta = \BF{a}\cdot\BF{b}_{\perp}$ and $\acute{\beta} = \BF{a}\cdot\BF{b}_{\perp}$, so $\beta = \acute{\beta}$.

	 	In the plane of the vectors $\BF{a}$ and $\BF{b}$, there are two vectors perpendicular to $\BF{b}$ (see Fig. \ref{fig1b} ), which we will call $\BF{b}_{\perp 1}$ and $\BF{b}_{\perp 2}$, where $\BF{b}_{\perp 1} = -\BF{b}_{\perp 2}$. 
	 	
	 	Since $\beta = \acute{\beta}$, we have that	 	
	 	 $\BF{a}\cdot \BF{b}_{\perp 1} = \BF{a}\cdot \BF{b}_{\perp 2} $. If we consider that $\BF{b}_{\perp 1}\neq \BF{b}_{\perp 2}$ we got to that
	 	$2\,\BF{a}\cdot\BF{b}_{\perp 1} = 0$, which is a contradiction because $\BF{a}\neq 0 $ and $\BF{b}_{\perp 1}$ is not perpendicular to $\BF{a}$ in the general sense. Therefore our assumption that $\BF{b}_{\perp 1} \neq \BF{b}_{\perp 2}$ is not correct, then, $\BF{b}_{\perp 1} = \BF{b}_{\perp 2}$.\\
	 	
	 	On the other hand we have that $R\,\BF{b}^{-1} = \BF{b}_{\perp 1}^{-1}$ and $\acute{R}\,\BF{b}^{-1} = \BF{b}_{\perp 2}^{-1} = \BF{b}_{\perp 1}^{-1}$ then we conclude that $R = \acute{R}$, so we show that $\BF{a}\otimes\BF{b} = \BF{b}\otimes\BF{a}$.

	 \end{proof}
\item[M3:]	
	 $\BF{a}\otimes\BF{u} = \BF{u}\otimes\BF{a}= a\,I$, identity element of multiplication.
	 \begin{proof}
	 	We are going to show that there exists a vector $\BF{u}$ such that $\BF{a}\otimes\BF{u}= a\,I$.\\
	 	
	 	By definition  $\BF{a}\otimes\BF{u} =\alpha\, I\ +\beta\, R$, where $\alpha =  \BF{a}\cdot \BF{u}$ and $\beta =  \BF{a}\cdot \BF{u}_{\perp} = \BF{a}_{\perp}\cdot \BF{u},$ but $\BF{a}\cdot \BF{u} = a$ and $\BF{a}\cdot \BF{u}_{\perp} = 0$,
	 	 so	$u = 1/cos \theta$ and $au\,tan\theta = 0$ then $\theta = 0$ and $u = 1$.\\
	 	 
	 	Since $\BF{u} = \left(\BF{a}\otimes\BF{u}\right)\BF{a}^{-1}  $
	 	 we have that
	 	 
	 	 $\BF{u} = a\,I\,\BF{a}^{-1}
	 	         = a\,\BF{a}/a^2 = \BF{a}/a$. 
	 	
	 	That is, we find that $\BF{u}= \BF{a}/a$.
	 \end{proof}
\item[M4:]
	 $\BF{a}\otimes(\BF{b}+\BF{c}) =\BF{a}\otimes\BF{b}+\BF{a}\otimes\BF{c} $, distributive with respect to vector addition.
	 
	 \begin{proof}
	 	By definition
	 	
	 	$\BF{a}\otimes(\BF{b}+\BF{c})= \BF{a}\cdot(\BF{b}+\BF{c})\,I + \BF{a}\cdot(\BF{b}+\BF{c})_{\perp}\,R$.
	 	
	 	Since the dot product is distributive, we have that
	 	
	 	$\BF{a}\otimes(\BF{b}+\BF{c})= \BF{a}\cdot\BF{b}\,I +\BF{a}\cdot\BF{c}\,I + \BF{a}\cdot\BF{b}_{\perp}\,R +\BF{a}\cdot\BF{c}_{\perp}\,R $\\
	 	
	 	$\BF{a}\otimes(\BF{b}+\BF{c})= (\BF{a}\cdot\BF{b}\,I + \BF{a}\cdot\BF{b}_{\perp}\,R )+ (\BF{a}\cdot\BF{c}\,I + \BF{a}\cdot\BF{c}_{\perp}\,R )$.\\
	 	
	 	Then
	 	
	 		$\BF{a}\otimes(\BF{b}+\BF{c})= \BF{a}\otimes\BF{b}+\BF{a}\otimes\BF{c}$.
	 	
	 \end{proof}
\item[M5:] $\BF{a}\otimes k\,\BF{b} = k\,\BF{a}\otimes\BF{b}$, linear with respect to scalar multiplication, where $k\in \mathbb{R}$.
\begin{proof}
	By definition, $\BF{a}\otimes k\,\BF{b} = \alpha\, I\ +\beta\, R $, where $\alpha = \BF{a}\cdot k\BF{b} = k \BF{a}\cdot \BF{b}$ and $\beta = k\BF{a}\cdot\BF{b}_{\perp}$ then
	
	 $\BF{a}\otimes k\,\BF{b} = k\BF{a}\cdot \BF{b}\, I\ + k \BF{a}\cdot \BF{b}_{\perp}\, R = k\left(\BF{a}\otimes\BF{b}\right)$
\end{proof}
\item[M6:] $\BF{a}\otimes\BF{b} = O$, then either $\BF{a} = \BF{0}$ or $\BF{b} = \BF{0}$. 
\begin{proof}
	$\BF{a}\otimes\BF{b} = O$, means that $\alpha = 0 $ and $\beta = 0$ or $ab\left(cos\theta - sin\theta\right) = 0$ for any $\theta$, then necessarily $a = 0$ or $b = 0$, which means that $\BF{a} = \BF{0}$ or $\BF{b} = \BF{0}$.
	
\end{proof}

\end{enumerate}

Let's see now two examples of multiplication.  For this we return to the vectors of the examples \ref{example1} and \ref{example2} .\\

\begin{exmp}\label{example3}
	Let \BF{a} = (3, -1) and \BF{b} = (2, 5) be two vectors in $\mathbb{R}^2$ find $\BF{a}\otimes\BF{b}$.\\
	
	$\BF{a}\otimes\BF{b} = \BF{a}\cdot\BF{b}\, I + \BF{a}\cdot\BF{b}_{\perp}\,R$, where $\BF{a}\cdot\BF{b} = 1$ and $\BF{a}\cdot\BF{b}_{\perp} = -17$ (we choose $\BF{b}_{\perp} = (-5, 2)$). 
	
	$R = \left[ \begin{array}{clr}
	0 &  -1 \\
	1
	& 0
	\end{array}\right]$, then 
	
	$\BF{a}\otimes\BF{b} = \left[ \begin{array}{clr}
	1 &  17 \\
	-17
	& 1
	\end{array}\right]$.
	
	If $\BF{b} = 1/29\, (2, 5)$ (the vector inverse of $\BF{b}$) then we can see that $\BF{a}\otimes\BF{b}$ match with the division of the example \ref{example1} as we hope.  
	
\end{exmp}

\begin{exmp}\label{example4}
	Let \BF{a} = (3, -1, 2) and \BF{b} = (2, 5, 1) be two vectors in $\mathbb{R}^3$ find $\BF{a}\otimes\BF{b}$.\\
	
	$\BF{a}\otimes\BF{b} = 3\,I + \sqrt{411}\,R$,
	
	were $\alpha = \BF{a}\cdot\BF{b} = 3$ and $\beta = \BF{a}\cdot\BF{b}_{\perp} = \sqrt{411}$
	
	\begin{eqnarray}
	R  &= &	\frac{1}{\sqrt{411}}\left[ \begin{array}{clrr}
	0 &  17 & -1\\
	-17 & 0 & -11\\
	1& 11& 0
	\end{array}\right],\nonumber
	\end{eqnarray}
	
	then
	
		\begin{eqnarray}
	\BF{a}\otimes\BF{b}  &= &	\left[ \begin{array}{clrr}
	3 &  17 & -1\\
	-17 & 3 & -11\\
	1& 11& 3
	\end{array}\right].\nonumber
	\end{eqnarray}
	
	If $\BF{b} = 1/30\,(2, 5, 1)$ (the inverse vector of $\BF{b}$), then we can see that $\BF{a}\otimes\BF{b}$ match with the division of the example \ref{example2} .
\end{exmp}

\section{Conclusions}

It is possible in the framework of vector algebra to define the division between vectors in the vector spaces $\mathbb{R}^2$ and $\mathbb{R}^3$, which is defined as a matrix. This matrix represents rotation together with scaling. It is possible to define the corresponding multiplication, which represents, in the same way, a rotation. The inverse vector was defined through this multiplication, through which it was possible to obtain the division. Fundamental properties of division and multiplication between vectors, analogous to the properties of operations with numbers, were given along with their proofs. We will continue with the study of these operations and their applications in Physics. We think that they can be linked to group theories.

\end{multicols*}

\end{document}